\documentclass[10pt,a4paper,leqno]{amsart}

\date{}

\usepackage{latexsym}

\usepackage{amsfonts}

\usepackage{amsmath}

\usepackage{amssymb}

\usepackage{amscd}

\usepackage{geometry}

\input xy
\xyoption{all}

\newcommand{\ma}[1]{\ensuremath{\mathbb{#1}}}
 \newcommand{\bo}[1]{\ensuremath{\bf{#1}}}
\font\bb=msbm7 at 10 pt

\def \C {\hbox{\bb C}}
\def \Z {\hbox{\bb Z}}
\def \Q {\hbox{\bb Q}}

\def \P {\mathcal{P}}
\def \F {\hbox{\bb F}}

\def \A {\hbox{\bb A}}

\def \Tr {\mbox{\rm{Tr}}}

\def \Vect {\mbox{\rm{Vect}}}

\def \T {\mathcal{T}}

\def \G {\mathcal{G}}

\def \L {\mathcal{L}}
\def \H {\mathcal{H}}
\def \U {\mathcal{U}}
\def \O {\mathcal{O}}
\def \K {\mathcal{K}}
\def \Q {\hbox{\bb Q}}
\def \I {\mbox{\bf I}}
\def \M {\mathcal{M}}

\def \sgn {\mbox{\rm{sgn}}}
\def \d {\mbox{\bf d}}

\newcommand{\Mat}{\ensuremath{\mbox{\rm{Mat}}}}

\begin{document}

\title{Newton polygons for twisted exponential sums and polynomials $P(x^d)$.}

 \keywords{$L$-functions of exponential sums, Newton
polygon, Newton strata, Hasse polynomials} 

\subjclass[2000]{11T23,11L03,14G15}

\author{R\'egis Blache}
\address{
\'Equipe ``G\'eom\'etrie Alg\'ebrique et Applications \`a la Th\'eorie de l'Information'', Universit\'e de Polyn\'esie Fran\c{c}aise, BP 6570, 98702 FAA'A, Tahiti, Polyn\'esie Fran\c{c}aise}
\email{blache@upf.pf}
\author{\'Eric F\'erard}
\address{
\'Equipe ``G\'eom\'etrie Alg\'ebrique et Applications \`a la Th\'eorie de l'Information'', Universit\'e de Polyn\'esie Fran\c{c}aise, BP 6570, 98702 FAA'A, Tahiti, Polyn\'esie Fran\c{c}aise}
\email{ferard@upf.pf}

\begin{abstract}
We study the $p$-adic absolute value of the roots of the $L$-functions associated to certain twisted character sums, and additive character sums associated to polynomials $P(x^d)$, when $P$ varies among the space of polynomial of fixed degree $e$ over a finite field of characteristic $p$. For sufficiently large $p$, we determine in both cases generic Newton polygons for these $L$-functions, which is a lower bound for the Newton polygons, and the set of polynomials of degree $e$ for which this generic polygon is attained. In the case of twisted sums, we show that the lower polygon defined in \cite{as1} is tight when $p\equiv 1~[de]$, and that it is the actual Newton polygon for any degree $e$ polynomial.
\end{abstract}

\maketitle

\section{Introduction}

Let $k:=\F_q$, $q:=p^m$ be a finite field, $P=X^e+a_{e-1}X^{e-1}+\dots+a_1X \in k[X]$ a degree $e$ polynomial over $k$; for any integer $r\geq 1$, we denote by $k_r=\F_{q^r}$ the extension of degree $r$ of $k$. We fix once and for all a non trivial additive character $\Psi$ of $\F_p$, with values in $\C_p^\times$; for any integer $r\geq 1$, the map $\Psi_{mr}:=\Psi\circ \Tr_{\ma{F}_{p^{mr}}/\ma{F}_p}$ is a non trivial additive character of $k_r$. On the other hand, let $d\geq 2$ be a divisor of $q-1$, and $\chi$ be a multiplicative character of order $d$ of $k^\times$; if $N_{k_r/k}$ denotes the norm from $k_r$ to $k$, the map $\chi\circ N_{k_r/k}$ is a multiplicative character of order $d$ of $k_r^\times$.

\medskip
 
Let $1\leq \kappa\leq d-1$ be an integer. To these data one associates the character sums (defined for any $r\geq 1$)

$$S_r(P,\chi^\kappa):=\sum_{x\in k_r^\times} \Psi_{mr}(P(x))\chi^\kappa\circ N_{k_r/k}(x)$$

and the $L$-function

$$L(P,\chi^\kappa;T):=\exp\left(\sum_{r\geq 1}\frac{S_r(P,\chi^\kappa)}{r}T^r\right).$$

It follows from the work of Weil on the Riemann hypothesis in characteristic $p$ that this function is a polynomial of degree $e$ ({\it cf.} for instance \cite{se}, \cite{wei}); moreover its reciprocal roots (in $\C$) are algebraic integers with complex absolute value $q^{\frac{1}{2}}$, and $\ell$-adic absolute value $1$ for any prime $\ell\neq p$. 

\medskip

We shall adress here the question of their $p$-adic absolute value; we study the Newton polygons associated to the $L$-functions above. These polygons encode the $p$-adic absolute values of the reciprocal roots in the following way: if the Newton polygon (with respect to the valuation $v_p$) has a segment of (horizontal) length $l$ and slope $s$, the $L$-function has exactly $l$ reciprocal roots $\alpha_1,\dots,\alpha_l$ with $v_p(\alpha_i)=s$.

\medskip

Much attention has been adressed to this question in the case of additive exponential sums ({\it cf.} \cite{wan} and references therein). In the case of twisted exponential sums, Adolphson and Sperber ({\it cf.} \cite{as1}, \cite{as2}) have studied the degree of the associated $L$-function when $f$ is a $n$-variable Laurent polynomial; they have also given a lower bound (\cite{as1} Theorem 3.20, \cite{as2} Corollary 3.18) for the Newton polygon associated to this $L$-function. In this paper, we refine these results in the case of a one-variable polynomial: we show that the bound in \cite{as1}, \cite{as2} is tight when $p\equiv 1~[de]$, and give the exact lower bound in the other cases.

\medskip

More precisely, we study the variation of the Newton polygon $NP_q(P,\chi^\kappa)$ of $L(P,\chi^\kappa;T)$ with respect to the $q$-adic valuation, when $P$ runs over polynomials of degree $e$ in $k[X]$. Let us recast our results in terms of crystals; let $\L_\Psi$ denote the {\it Artin Schreier crystal} over $\A^1$, and $\H_{\chi^\kappa}$ denote the {\it Kummer crystal} over $\A^1\backslash\{0\}$ associated to the character $\chi^\kappa$; these are overconvergent $F$-isocrystals ({\it cf.} \cite{els} 6.5), and for any polynomial $P\in k[x]$ of degree $e$, we have an overconvergent $F$-isocrystal $P^*\L_\Psi\otimes\H_{\chi^\kappa}$ with 
$$L(P,\chi^\kappa;T)=\det\left(1-TF^*|H^1_{\rm rig,c}(\A^1\backslash\{0\}/K,P^*\L_\Psi\otimes\H_{\chi^\kappa})\right).$$
If we parametrize the set of degree $e$ monic polynomials without constant coefficient by the affine space $\A^{e-1}$, associating the point $(a_1,\dots,a_{e-1})$ to the polynomial $P(X)=X^e+a_{e-1}X^{e-1}+\dots+a_1X$, we can consider the family of overconvergent $F$-isocrystals $P^*\L_\Psi\otimes\H_{\chi^\kappa}$. Grothendieck's specialization theorem ({\it cf.} \cite{gr}, \cite{ka} Corollary 2.3.2) asserts that there is a Zariski dense open subset $U_{d,e,q,\kappa}$, the {\it open Newton stratum}, of the affine space $\A^{e-1}$, and a {\it generic Newton polygon} $GNP(d,e,q,\kappa)$ such that for any $P\in  U_{d,e,q,\kappa}$, we have $GNP(d,e,q,\kappa)=NP_q(P,\chi^\kappa)$. Moreover for any monic polynomial of degree $e$, $NP_q(P,\chi^\kappa)$ lies above the generic Newton polygon.

\medskip

We shall determine exactly the generic Newton polygon; this refines the lower bound in \cite{as1}, \cite{as2}. We also give the {\it Hasse polynomial} $P_{d,e,q,\kappa}$, i.e. the polynomial such that $U_{d,e,q,\kappa}$ is the complementary of the hypersurface $P_{d,e,q,\kappa}=0$. We show that these data only depend on $p$, i.e. are independent of the choice of the power $q$. Similar questions have been adressed in \cite{bf}, in the case of purely additive exponential sums associated to polynomials.

\medskip

Considering Newton polygons for exponential sums often leads to consider a combinatorial polygon (the {\it Hodge polygon}) which is a lower bound for the Newton polygons; in the present situation, we also get such a polygon, but since its slopes have something to do with the valuations of Gauss sums, we shall call it th {\it Hodge-Stickelberger polygon}. We now describe it: let $1\leq r\leq d-1$ be an integer prime to $d$; multiplication by $r$ modulo $d$ induces a permutation of the set $\{1,\dots,d\}$, i.e. an element $\sigma_r$ of $S_d$, the $d$-th symmetric group. For any $1\leq \kappa\leq d-1$, let $<\kappa>_r$ be the cycle of $\sigma_r$ containing $\kappa$. We define $$\mu_\kappa:=\frac{\sum_{i\in <\kappa>_r} i}{d\#<\kappa>_r};$$ 
this is a rational number in the interval $]0,1[$; note that it appears in Stickelberger's theorem on the valuation of Gauss sums ({\it cf.} Remark 3.1). We define the {\it Hodge-Stickelberger polygon} $HS(d,e,r,\kappa)$ as the joint of segments of (horizontal) length $1$, and slopes $\frac{i+\mu_{d-\kappa}}{e}$ for $0\leq i\leq e-1$. 

\medskip

This polygon has already been defined in an other way ({\it cf.} Remark 3.1) in \cite{as1} Theorem 3.20, \cite{as2} Corollary 3.18, where it is shown that the Newton polygon of the $L$-function with respect to the $q$-adic valuation always lies on or above it. Here we find this result in an other way, showing that if $r$ is the least positive residue of $p$ modulo $d$, then the generic Newton polygon lies on or above the Hodge-Stickelberger polygon. Moreover, we show the equality when $p\equiv 1~[de]$; in this case the situation is particularly simple, since for any polynomial of degree $e$, the Newton polygon of the associated $L$-function is exactly the Hodge-Stickelberger polygon.

\medskip

The situation is very similar to the case of $L$-functions associated to additive character sums, and polynomials of degree $e$: recall that the Newton polygons of these functions lie above the Hodge polygon $H(e)$ of length $e-1$ and slopes $\frac{i}{e}$, $1\leq i\leq e-1$, and that they are equal when $p\equiv 1~[d]$ ({\it cf.} \cite{ro}). Thus the Hodge-Stickelberger polygon $HS(d,e,r,\kappa)$ is obtained from the classical one by adding a segment of slope $\frac{\mu_{d-\kappa}}{e}$, and $\frac{\mu_{d-\kappa}}{e}$ to each of the slopes.

\medskip

Second, we shall look at Newton polygons of $L$-functions associated to an additive character and a polynomial of the form $P(x^d)$; here we do not need no more the hypothesis that $d$ divides $q-1$. If $\Psi$ and $P$ are as above, we define for any $r\geq 1$

$$S_r(P(x^d)):=\sum_{x\in k_r} \Psi_{mr}(P(x^d))$$

and the $L$-function

$$L(P(x^d);T):=\exp\left(\sum_{r\geq 1}\frac{S_r(P(x^d))}{r}T^r\right).$$

These $L$-functions have a natural link with the $L$-functions of twisted sums studied above ({\it cf.} Proposition 4.1). In particular their Newton polygons $NP(P(x^d))$ are concatenations of Newton polygons as above, and the same is true for the generic Newton polygon $GNP(d,e,p)$, and the Hodge-Stickelberger polygon $HS(d,e,r)$ associated to this situation. In the case $p\equiv 1~[de]$, our polygon coincides with the classical Hodge one $H(de)$, and we find a particular case of a theorem of Robba ({\it cf.} \cite{ro} Th\'eor\`eme 7.5): for any polynomial $P(x^d)$, the Hodge-Stickelberger polygon coincides with the Newton polygon $NP(P(x^d))$. In the other cases, the situation is much more intricated, and we determine both the generic Newton polygon and the Hasse polynomial for sufficiently great $p$.

\medskip

We say a few words about the asymptotic behaviour of these Newton polygons; it is not the aim of this article, but it has drown much attention recently ({\it cf.} \cite{wan}, Conjectures 1.11 and 1.12), and some material from the present article can be used to derive results in this fashion. Moreover this sheds some light on what happens in the non generic case. It is known that for $L$-functions associated to one variable additive exponential sums, the generic Newton polygon tends to the Hodge polygon when $p$ tends to infinity. More precisely ({\it cf.} \cite{zhu1}, \cite{zhu2}), for any $P\in \overline{\Q}[X]$, let the polygon $NP_q(P)$ be the Newton polygon obtained from the reduction of $P$ modulo a prime above $p$ in the field defined by the coefficients of $P$. Then there is a Zariski dense open subset $\U$ defined over $\Q$ in the space of polynomials of degree $e$ such that, for every $P$ in $\U(\overline{\Q})$, we have $\lim_{p\rightarrow \infty} NP_q(P)=H(d)$. 

\medskip

One cannot hope such a result here: since there are many Hodge-Stickelberger polygons (one for each residue prime to $d$ in $\Z/d\Z$) there is no such limit when $p$ tends to infinity. But we have the following weaker result, that we shall not prove here: when $p$ tends to infinity, and satisfies the additional condition $p\equiv r~[d]$, the generic Newton polygon $GNP(d,e,p)$ tends to the Hodge-Stickelberger polygon $HS(d,e,r)$ . Actually there is a Zariski dense open subset $\U_{d,e,r}$ defined over $\Q$ in the space of polynomials of degree $e$ over $\overline{\Q}$ such that for any $P$ in $\U_{d,e,r}(\overline{\Q})$ the Newton polygon $NP(P(x^d))$ tends to the Hodge-Stickelberger polygon $HS(d,e,r)$ when $p$ tends to infinity along the class of $r$ modulo $d$.

\medskip

Let us describe the methods we employ. We use $p$-adic cohomology, following the ideas of Dwork; actually our tools are mainly inspired by the work of Robba ({\it cf.} \cite{ro}). To be more precise, we use Washnitzer-Monsky spaces of overconvergent series $\H^\dagger(A)$, where $A$ is a $p$-adic annulus; one can define linear operators $\beta^{(0)},\dots,\beta^{(m-1)}$ on $\H^\dagger(A)$ and differential operators $D^{(0)},\dots,D^{(m-1)}$ with finite index on this space such that the $\beta^{(i)}$ and $D^{(i)}$ commute up to $p$. Then the linear map $\alpha=\beta^{(m-1)}\circ\beta^{(m-2)}\circ\dots\circ\beta^{(0)}$ commutes with $D^{(0)}$ up to $q$, and induces a linear map $\overline{\alpha}$ on the quotient $\H^\dagger(A)/D^{(0)}\H^\dagger(A)$ which has characteristic polynomial equal to $L(P,\chi^\kappa;T)$. Using monomial bases for the spaces $\H^\dagger(A)/D^{(s)}\H^\dagger(A)$, we are able to give congruences for the coefficients of the matrices $M:=\Mat(\overline{\beta}^{(s)})$ in terms of the coefficients of a lift of $P$. We deduce congruences for the minors of $N:=\Mat(\overline{\alpha})$, i.e. for the coefficients of the function $L(P,\chi^\kappa;T)$.

\medskip

The article is organized as follows: in section 1, we present the results from $p$-adic cohomology that we shall use throughout the paper. We define an operator $\overline{\alpha}$ over a finite dimensional $p$-adic vector space, whose characteristic polynomial is the $L$-function, and we describe a factorisation of $\overline{\alpha}$ in terms of linear maps whose matrices with respect to certain bases have a simple description. Then in section 2 we give $p$-adic approximations for the coefficients of these matrices, and for their minors, which in turn give the ``principal" $p$-adic parts of the coefficients of the $L$-function; we also express the Hasse polynomial, in order to determine the open Newton stratum. In section 3, we describe the generic and Hodge-Stickelberger polygons associated to twisted character sums. Finally, in section 4, we give the link between the $L$-function $L(P(x^d);T)$ and the $L$-functions $L(P,\chi^\kappa;T)$. We deduce results about the Newton polygons $NP(P(x^d))$, and describe the generic and Hodge-Stickelberger polygons attached to this situation.

\medskip

\section{$p$-adic theory.} 

The aim of this section is to express the $L$-function $L(P,\chi^\kappa;T)$ as the characteristic polynomial of a linear map acting on a $p$-adic vector space of dimension $e$. We shall only describe the pieces from Dwork's theory and $p$-adic cohomology we need here. The reader interested in more details can refer to \cite{dw}, \cite{ro} for instance.

\medskip

\subsection{p-adic analytic theory.} We begin by recalling some elements from Dwork's theory. Let $\Q_p$ be the field of $p$-adic numbers; we fix once and for all an algebraic closure $\overline{\Q}_p$ and a completion $\C_p$ of this algebraic closure; note that the field $\C_p$ remains algebraically closed. For any integer $m$, let $\K_m$ denote the unramified extension of degree $m$ of $\Q_p$ in $\C_p$. There are exactly $q-1$ elements of finite order in $\K_m^\times$; they form a multiplicative group, and we shall denote it by $\T_m^\times$. The set $\T_m:=\T_m^\times\cup\{0\}$ is called the {\it Teichm\"uller} of $\K_m$. Let $\O_m$ denote the valuation ring of $\K_m$; this is a local ring with maximal ideal generated by $p$, and residue field isomorphic to $k$. Reduction modulo $p$ induces a bijection from $\T_m$ to $k$, and we denote by $a\mapsto \widehat{a}$ its inverse; this map is called {\it Teichm\"uller lifting}. 

\medskip

For any $x\in \C_p$, $r>0$, let $B(x,r^+)$ ({\it resp.} $B(x,r^-)$) denote the closed ({\it resp.} open) ball in $\C_p$ with center $x$ and radius $r$; let $A$ be the annulus $B(0,1^+)\backslash B(0,1^-)$ in $\C_p$. We consider the space $\H^\dagger(A)$ of overconvergent analytic functions on $A$. It is well known since Dwork that one can represent analytically the characters of $\F_p$ and its extensions with elements in $\H^\dagger(A)$. Since $\Psi$ is a non trivial additive character of $\F_p$, $\zeta_p:=\Psi(1)$ is a non trivial $p$-th root of unity. There is an unique root $\pi$ of $X^{p-1}+p$ in $\C_p$ such that $\zeta_p-1\equiv \pi~ [\pi^2]$. The power series $\theta(X)$ coming from the function $\exp(\pi X-\pi X^p)$ has radius of convergence $>1$, and the map from $\F_p$ to $\C_p^\times$ defined by $x\mapsto \theta(\widehat{x})$ is exactly the character $\Psi$. Moreover, we define for any $m\geq 1$ the power series $\theta_m(x):=\theta(x)\theta(x^p)\dots\theta(x^{p^{m-1}})$ (note that this is the power series associated to the function $\exp(\pi X-\pi X^{q})$); once again this is an overconvergent power series, and the map $x\mapsto \theta_m(\widehat{x})$ is exactly the character $\Psi_m$. These functions first appeared in Dwork \cite{dw}, and are often called {\it splitting functions}. On the other hand, for any $r\geq 1$ $\T_{mr}^\times$ is the group of $q^r-1$-th roots of unity in $\C_p^\times$, and the map $x\mapsto \widehat{x}$ from $k_r$ to $\T_{mr}\subset \K_{mr}^\times$ is a multiplicative character of order $q^{r}-1$; consequently the map $x\mapsto \widehat{x}^{\frac{q^{r}-1}{d}}$ is a multiplicative character of order $d$ of $k_r$, and we shall denote it by $\chi_d^{(r)}$. Since Gal$(\K_{mr}/\K_m)$ is a cyclic group, generated by the Frobenius acting as the $q$-th power on $\T_{mr}$, we see that the norm from $\K_{mr}$ to $\K_m$ acts as the $(q^r-1)/(q-1)$-th power on $\T_{mr}$, and the character $\chi_d^{(r)}$ equals $\chi_d \circ N_{k_r/k}$, where $\chi_d:=\chi_d^{(1)}$.

\medskip

We define the function $H(X):=X^{\frac{(q-1)\kappa}{d}}\prod_{i=1}^e \theta_m(\widehat{a_i}X^i)$. From the discussion above, this is an element of $\H^\dagger(A)$, and for any $x\in k^\times$, we have $H(\widehat{x})=\chi_d(x^\kappa)\Psi_m(P(x))$. Recall that Dwork's operator $\psi_q$ is the endomorphism of $\H^\dagger(A)$ defined by $\psi_qf(x):=\frac{1}{q}\sum_{z^q=x}f(z)$ (note that when $f(X)=\sum b_nX^n$, then $\psi_q f(X)=\sum b_{qn}X^n$). Now let $\alpha$ be the endomorphism $\psi_q\circ H$ of $\H^\dagger(A)$; by Monsky-Reich trace formula we have, for any $r\geq 1$
$$S_r(P,\chi_d^\kappa)=(q^{r}-1)\Tr(\alpha^r),$$
and we deduce the following expression for the $L$-function 
$$L(P,\chi_d^\kappa;T)=\frac{\det(1-T\alpha)}{\det(1-qT\alpha)}.$$

\subsection{$p$-adic cohomology} Let $\P(X):=X^e+\sum_{i=1}^{e-1} \widehat{a_i}X^i$ be the polynomial in $\K_m[X]$ obtained by Teichm\"uller lifting of the coefficients of $P$. Consider the differential operator $D:=X\frac{d}{dX}+\pi X\P'-\frac{\kappa}{d}$ acting on $\H^\dagger(A)$; it is injective since a solution of $Df=0$ is the function $F(X)=X^{\frac{\kappa}{d}}\exp(-\pi \P(X))$ which is not overconvergent. Moreover its cokernel has dimension $e$ from \cite{ro} Th\'eor\`eme 6.1, and we have the commutative diagram with exact rows
$$\xymatrix{
 0 \ar[r]& \ar[d]_{q\alpha} \H^\dagger(A) \ar[r]^{D} & \ar[d]_{\alpha} 
\H^\dagger(A) \ar[r] & \ar[d]_{\overline{\alpha}} \H^\dagger(A)/D\H^\dagger(A) \ar[r] & 0\\
0 \ar[r]& \H^\dagger(A) \ar[r]^{D} & \H^\dagger(A) \ar[r] & \H^\dagger(A)/D\H^\dagger(A) \ar[r] & 0,\\
}$$
where $\overline{\alpha}$ is the morphism of $\H^\dagger(A)/D\H^\dagger(A)$ obtained from $\alpha$. From the discussion above we get $L(P,\chi_d^\kappa,T)=\det(1-T\overline{\alpha})$.

\medskip

Now we factorize the morphisms $\alpha$ and $\overline{\alpha}$, in order to get expressions for the matrix of $\overline{\alpha}$ with respect to a suitable basis of $\H^\dagger(A)/D\H^\dagger(A)$. For any $0\leq s\leq m$, we define $\kappa_{s}$ as the least positive integer such that $p^s\kappa_{s}\equiv \kappa~[d]$, and 
$$K_s:=\frac{p\kappa_{s+1}-\kappa_s}{d}.$$
Note that $\kappa_m=\kappa_0=\kappa$ since $q\equiv 1~[d]$. Let $\tau$ be a lifting of the Frobenius of $\overline{\F}_p$ to $\C_p$ such that $\tau(\pi)=\pi$. We define, for any $0\leq s\leq m-1$, the power series $H^{(s)}(X):=X^{K_s}\prod_{i=1}^e\theta(\widehat{a_i}^{\tau^s}X^i)$; this is an overconvergent power series, and we denote by $\beta^{(s)}$ the endomorphism $\psi_p\circ H^{(s)}$ of $\H^\dagger(A)$. 

\medskip

{\bf Lemma 1.1.} {\it We have the following factorisation for the morphism $\alpha$  $$\alpha=\beta^{(m-1)}\circ\dots\circ \beta^{(0)}.$$}

\medskip

{\it Proof.} Recall that for any analytic overconvergent function $f(X)\in \H^\dagger(A)$, we have the following commutation rule with Dwork's operator $f(X)\circ \psi_p =\psi_p\circ f(X^p)$; we can rewrite the operator $\beta^{(m-1)}\circ\dots\circ \beta^{(0)}$ in the following way
$$\begin{array}{rcl}
\beta^{(m-1)}\circ\dots\circ \beta^{(0)} & = & \psi_p\circ H^{(m-1)}(X)\circ\dots\circ  \psi_p\circ H^{(0)}(X)\\
& = & \psi_p\circ\dots\circ\psi_p \circ\left(H^{(m-1)}(X^{p^{m-1}})\dots H^{(0)}(X)\right)\\
& = & \psi_q \circ\left(H^{(m-1)}(X^{p^{m-1}})\dots H^{(0)}(X)\right).\\
\end{array}$$
Now we have 
$$H^{(m-1)}(X^{p^{m-1}})\cdots H^{(0)}(X)=X^{K_0+\dots+p^{m-1}K_{m-1}}\prod_{s=0}^{m-1}\prod_{i=1}^e\theta(\widehat{a_i}^{\tau^s}X^{p^s i});$$
from the definition of the $K_s$, $K_0+\dots+p^{m-1}K_{m-1}=\frac{(q-1)\kappa}{d}$, and the second product is exactly $\prod_{i=1}^e \theta_m(\widehat{a_i}X^i)$. Thus we get $H^{(m-1)}(X^{p^{m-1}})\cdots H^{(0)}(X)=H(X)$, and this ends the proof of the lemma.

\medskip

Let $D^{(s)}$ be the differential operator $X\frac{d}{dX}+\pi X\P'^{\tau^s}-\frac{\kappa_s}{d}$ acting on $\H^\dagger(A)$; once again this is an injective morphism, and its cokernel has dimension $e$ (note that $D^{(m)}=D^{(0)}=D$). For any $0\leq s\leq m-1$ we have the commutative diagram 
$$\xymatrix{
 0 \ar[r]& \ar[d]_{p\beta^{(s)}} \H^\dagger(A) \ar[r]^{D^{(s)}} & \ar[d]_{\beta^{(s)}} 
\H^\dagger(A) \ar[r] & \ar[d]_{\overline{\beta^{(s)}}} \H^\dagger(A)/D^{(s)}\H^\dagger(A) \ar[r] & 0\\
0 \ar[r]& \H^\dagger(A) \ar[r]^{D^{(s+1)}} & \H^\dagger(A) \ar[r] & \H^\dagger(A)/D^{(s+1)}\H^\dagger(A) \ar[r] & 0,\\
}$$
and we obtain the following factorisation for $\overline{\alpha}$: $\overline{\alpha}=\overline{\beta^{(m-1)}}\circ\cdots\circ\overline{\beta^{(0)}}$.

\medskip

\section{$p$-adic approximation of the coefficients of the $L$-function.} 

The aim of this section is to give congruences for the coefficients of the matrices defined in the former section, and to deduce congruences for the coefficients of the $L$-functions $L(P,\chi^\kappa;T)$.

\medskip

\subsection{The matrices $B^{(s)}$ and their coefficients} We now describe the matrices of the endomorphisms $\beta^{(s)}$. We begin by describing bases for the spaces $\H^\dagger(A)/D^{(s)}\H^\dagger(A)$. 

\medskip

First remark that the differential operators $D^{(s)}$ have their coefficients in $\K_m(\pi)[X,X^{-1}]$; from the comparison theorem between algebraic and analytic cohomology ({\it cf.} \cite{ro}), we see that a complementary subspace for $D^{(s)} \K_m(\pi)[X,X^{-1}]$ in $\K_m(\pi)[X,X^{-1}]$ is a complementary subspace for $D^{(s)}\H^\dagger(A)$ in $\H^\dagger(A)$. Now for any $n\in \Z$ we have
$$ D^{(s)}(X^{n-e})=(n-e-\frac{\kappa_s}{d})X^{n-e}+\pi\sum_{i=1}^e i\widehat{a_i}^{\tau^s}X^{n-e+i},$$
and the monomials $X,\dots,X^e$ span a complementary subspace for $D^{(s)}\H^\dagger(A)$ in $\H^\dagger(A)$. For any $n\in \Z$, $0\leq s\leq m$, we set $X^n\equiv\sum_{i=1}^e a_{n,i}^{(s)}X^i~[D^{(s)}\H^\dagger(A)]$. On the other hand set $G(X)=\prod_{i=1}^e\theta(\widehat{a_i}X^i):=\sum_{n\geq 0} h_nX^n$; since $\theta(X)\in \Q_p(\zeta_p)[[X]]$, we have $\prod_{i=1}^e\theta(\widehat{a_i}^{\tau^s}X^i)=G(X)^{\tau^{s}}=\sum_{n\geq 0} h_n^{\tau^{s}}X^n$, and $H^{(s)}(X)=\sum_{n\geq K_s} h_{n-K_s}^{\tau^{s}}X^n$, where we set $h_n:=0$ for negative $n$. Thus we get, for any $1\leq j\leq e$, $$\beta^{(s)}(X^j)=\psi_p(H^{(s)}(X)X^j)=\psi_p\left(\sum_{n\geq K_s+j} h_{n-K_s-j}^{\tau^{s}}X^n\right)=\sum_{n\geq l_s} h_{pn-K_s-j}^{\tau^{s}}X^n,$$
with $l_s=\lceil\frac{K_s+j}{p}\rceil$, and passing to the cohomology (in $\H^\dagger(A)/D^{(s+1)}\H^\dagger(A)$)
$$\begin{array}{rcl}
\beta^{(s)}(X^j)& \equiv & \sum_{n\geq l_s} h_{pn-K_s-j}^{\tau^{s}}\sum_{i=1}^e a_{n,i}^{(s+1)}X^i~[D^{(s+1)}\H^\dagger(A)]\\
& \equiv & \sum_{i=1}^e \sum_{n\geq l_s} h_{pn-K_s-j}^{\tau^{s}}a_{n,i}^{(s+1)}X^i~[D^{(s+1)}\H^\dagger(A)]\\
& \equiv & \sum_{i=1}^e \left( h_{pi-K_s-j}^{\tau^{s}}+\sum_{n\geq e+1} h_{pn-K_s-j}^{\tau^{s}}a_{n,i}^{(s+1)}\right)X^i~[D^{(s+1)}\H^\dagger(A)],\\
\end{array}$$
the last congruence being a consequence of the equality $a_{n,i}^{(s+1)}=\delta_{ni}$ for $1\leq n\leq e$. Thus the matrix $B^{(s)}$ of $\beta^{(s)}$ with respect to the bases $\{X,\dots,X^e\}$ of $\H^\dagger(A)/D^{(s)}\H^\dagger(A)$ and $\H^\dagger(A)/D^{(s+1)}\H^\dagger(A)$ respectively can be written $(b_{ij}^{(s)})_{1\leq i,j\leq e}$, where 
$$b_{ij}^{(s)}=h_{pi-K_s-j}^{\tau^{s}}+\sum_{n\geq e+1} h_{pn-K_s-j}^{\tau^{s}}a_{n,i}^{(s+1)}.$$

\medskip

We  now give congruences for the coefficients of the matrices $B^{(s)}$; these will help us in the next section to determine the ``principal parts" of the coefficients of the $L$-function $L(P,\chi_d^\kappa;T)$. We begin by recalling some results from \cite{bf} on the coefficients $h_n$. Recall that $\P$ is the polynomial obtained by lifting the coefficients of $P$ to the Teichm\"uller of $\K_m$; in the following, for $f$ a polynomial, we denote by $\{f\}_n$ the coefficient of degree $n$ of $f$

\medskip

{\bf Proposition 2.1.} {\it Set $\nu^{(s)}_{i,j}:=\lceil \frac{pi-K_s-j}{e}\rceil$; then for any $1\leq i,j\leq e$, $0\leq s\leq m-1$, and any $p\geq 2de$ we have the congruence
$$b_{ij}^{(s)}\equiv \left\{\P^{\nu^{(s)}_{i,j}}\right\}_{pi-K_s-j}^{\tau^s}\frac{\pi^{\nu^{(s)}_{i,j}}}{\nu^{(s)}_{i,j}!}~\left[\pi^{\nu^{(s)}_{i,j}+1}\right].$$}

\medskip

{\it Proof.} First remark from the definition of the $\kappa_s$ that for any $0\leq s\leq m$ we have $1\leq \kappa_s\leq d-1$; thus $\frac{p-d+1}{d}\leq K_s\leq \frac{p(d-1)-1}{d}$, and we get $\frac{p+1}{d}-e\leq pi-K_s-j \leq p(e-\frac{1}{d})-\frac{1}{d}-2$ for any $1\leq i,j\leq e$. From the hypothesis $p\geq 2de$, we deduce the inequalities $1\leq pi-K_s-j \leq (p-1)e$ and $\nu^{(s)}_{i,j}\leq p-1$. Thus we can apply Lemma 2.1 in \cite{bf} to $h_{pi-K_s-j}^{\tau^{s}}$, and we get 
$$h_{pi-K_s-j}^{\tau^{s}}\equiv \left\{\P^{\nu^{(s)}_{i,j}}\right\}_{pi-K_s-j}^{\tau^s}\frac{\pi^{\nu^{(s)}_{i,j}}}{\nu^{(s)}_{i,j}!}~\left[\pi^{\nu^{(s)}_{i,j}+1}\right].$$
From the expression of $b_{ij}^{(s)}$ above, it remains to show that for any $n\geq e+1$ we have $h_{pn-K_s-j}^{\tau^{s}}a_{n,i}^{(s+1)}\equiv 0~[\pi^{\nu^{(s)}_{i,j}+1}]$. First remark that exactly as in the case $\kappa=0$ treated in \cite{ro} Lemma on p241, we have $v_\pi(a_{n,i}^{(s+1)})\geq-\left[\frac{n-i}{e}\right]$.

\medskip

Assume first $e+1\leq n\leq p$; since $pn-K_s-j<p^2$, we have $v_\pi(h_{pn-K_s-j})\geq \frac{pn-K_s-j}{e}$. If $n<e+i$,  $v_\pi(a_{n,i}^{(s+1)})\geq 0$, and $v_\pi(h_{pn-K_s-j}^{\tau^{s}}a_{n,i}^{(s+1)})\geq \frac{pn-K_s-j}{e}\geq\frac{pi-K_s-j}{e}+\frac{p}{e}$; now for $e+i\leq n\leq p$, $v_\pi(h_{pn-K_s-j}^{\tau^{s}}a_{n,i}^{(s+1)})\geq \frac{pn-K_s-j}{e}-\left[\frac{n-i}{e}\right]\geq\frac{(p-1)n-K_s+i-j}{e}\geq \frac{pi-K_s-j}{e}+p-1$, and in both cases we get $v_\pi(h_{pn-K_s-j}^{\tau^{s}}a_{n,i}^{(s+1)})\geq \nu^{(s)}_{i,j}+1$.

\medskip

For $n>p$, we have $v_\pi(h_{pn-K_s-j})\geq \frac{pn-K_s-j}{e}\left(\frac{p-1}{p}\right)^2$, and 
$$\begin{array}{rcl}
v_\pi(h_{pn-K_s-j}^{\tau^{s}}a_{n,i}^{(s+1)})& \geq & \frac{pn-K_s-j}{e}\left(\frac{p-1}{p}\right)^2-\frac{n-i}{e}\\
 & \geq & \frac{n}{e}\left(\frac{(p-1)^2}{p}-1\right)-\frac{1}{e}\left(\left(\frac{p-1}{p}\right)^2(K_s+j)-i\right).\\
 \end{array}$$
From the majoration of $K_s$ above, and since $\left(\frac{p-1}{p}\right)^2\leq 1$, we get $\left(\frac{p-1}{p}\right)^2(K_s+j)-i\leq e+p-\frac{p+1}{d}\leq p$ since $p\geq de$. Now since $n>p$, we have $v_\pi(h_{pn-K_s-j}^{\tau^{s}}a_{n,i}^{(s+1)}) \geq \frac{p}{e}\left(\frac{(p-1)^2}{p}-2\right)$, and this number is strictly greater than $p-1$ as soon as $p\geq e+4$; since $p\geq 2de$, this ends the proof of the proposition.

\medskip

\subsection{The coefficients of the $L$-function} We come to the sum of the 
$n\times n$ minors centered on the diagonal of $B^{(s)}$
$$M_n^{(s)}=\sum_{1\leq u_1<\dots<u_n\leq d-1} \sum_{\sigma\in S_n} 
\sgn(\sigma) \prod_{i=1}^n b^{(s)}_{u_iu_{\sigma(i)}}.$$

This is the coefficient of degree $n$ of the characteristic polynomial $\det(1-T\beta^{(s)})$; the evaluation of these minors will in turn allow us to give congruences for the coefficients of the $L$-function $L(P,\chi_d^\kappa;T)$.

\medskip

{\bf Definition 2.1.} {\it {\it i)} For any $1\leq n\leq e$, set $Y_n^{(s)}(\kappa):= \min_{\sigma\in S_n} 
\sum_{k=1}^n \nu^{(s)}_{k,\sigma(k)}$, and
$$\Sigma_n^{(s)}(\kappa):=\{\sigma\in S_n,~\sum_{k=1}^n 
\nu^{(s)}_{k,\sigma(k)}=Y_n^{(s)}(\kappa)\}.$$

{\it ii)} For every $1\leq i \leq d-1$, set $1\leq j_i^{(s)}(\kappa)\leq e$ be the least positive 
integer congruent to $pi-K_s$ modulo $e$, and for every $1\leq n \leq e$, 
let $B_n^{(s)}(\kappa):=\{1\leq i\leq n,~j_i^{(s)}(\kappa)\leq n\}$.}

\medskip

Note that since $p$ is coprime to $e$, the map $i\mapsto j_i^{(s)}(\kappa)$ is an 
element of $S_{e}$, the $e$-th symmetric group.  We can use the set 
$B_n^{(s)}(\kappa)$ to describe $\Sigma_n^{(s)}(\kappa)$ and evaluate $Y_n^{(s)}(\kappa)$ precisely.

\medskip

{\bf Lemma 2.2.} {\it Let $1\leq n\leq e$; we have $\Sigma_n^{(s)}(\kappa)=\{ 
\sigma\in S_n,~\forall i\in B_n^{(s)}(\kappa)~\sigma(i)\geq j_i^{(s)}(\kappa)\}$, and 
$Y_n^{(s)}(\kappa)=\sum_{k=1}^n \lceil\frac{pk-K_s}{e}\rceil-\#B_n^{(s)}(\kappa)$.}

\medskip

{\it Proof.} It is easily seen that for any $1\leq j\leq j_i^{(s)}(\kappa)-1$, we have 
$\nu^{(s)}_{i,j}=\lceil\frac{pi-K_s}{e}\rceil$, and for $j_i^{(s)}(\kappa)\leq 
j\leq e$, $\nu^{(s)}_{i,j}=\lceil\frac{pi-K_s}{e}\rceil-1$. From 
this we deduce, for any $\sigma\in S_n$
$$\sum_{k=1}^n \nu^{(s)}_{k,\sigma(k)}=\sum_{k=1}^n 
\lceil\frac{pk-K_s}{e}\rceil-\#\left\{1\leq k\leq n,~\sigma(k)\geq j_k^{(s)}(\kappa)\right\}.$$
Now we have the inclusion $\{1\leq k\leq n,~\sigma(k)\geq j_k^{(s)}(\kappa)\}\subset 
B_n^{(s)}(\kappa)$. Finally the set $\{ \sigma\in S_n,~\forall i\in B_n^{(s)}(\kappa)~\sigma(i)\geq j_i^{(s)}(\kappa)\}$ is not empty, since $i\mapsto j_i^{(s)}(\kappa)$ is an injection from $B_n^{(s)}(\kappa)$ 
into $\{1,\dots,n\}$; we get $Y_n^{(s)}(\kappa)=\sum_{k=1}^n 
\lceil\frac{pk-K_s}{e}\rceil-\#B_n^{(s)}(\kappa)$, and the permutations reaching this 
minimum  are exactly the ones with $\sigma(i)\geq j_i^{(s)}(\kappa)$ for all $i\in 
B_n^{(s)}(\kappa)$. This is the desired result.

\medskip

{\bf Definition 2.2.} {\it Recall that we have set $\P(X)=\sum_{i=1}^e 
\widehat{a_i}X^i$. For any $1\leq n\leq e$ let $\P_{n,\kappa}^{(s)}$ be the polynomial in 
$\Z[X_1,\dots,X_e]$ defined by
$$\P_{n,\kappa}^{(s)}(\widehat{a_1},\dots,\widehat{a_e}):=\sum_{\sigma\in \Sigma_n^{(s)}(\kappa)} \sgn(\sigma)
\prod_{i=1}^n\left\{\P^{\nu^{(s)}_{i,\sigma(i)}}\right\}_{pi-\sigma(i)-K_s}.$$}

\medskip

We are now ready to give a congruence for the coefficients $M_n^{(s)}$ of the 
polynomial $\det(1-T\beta^{(s)})$. We don't rewrite the proof since it is very similar to the one of Proposition 2.2 in \cite{bf}.

\medskip

{\bf Proposition 2.2.} {\it Assume $p\geq 2de$; then for any $1\leq n\leq 
e$, we have
$$M_n^{(s)}\equiv \frac{\P_{n,\kappa}^{(s)}(\widehat{a_1},\dots,\widehat{a_e})^{\tau^s}}{\prod_{i\notin 
B_n^{(s)}(\kappa)}\lceil\frac{pi-K_s}{e}\rceil!\prod_{i\in 
B_n^{(s)}(\kappa)}\left(\lceil\frac{pi-K_s}{e}\rceil-1\right)!}\pi^{Y_n^{(s)}(\kappa)}\quad[\pi^{Y_n^{(s)}(\kappa)+1}].$$}

\medskip

Finally, using the same methods as in Section 3 in \cite{bf}, we get a congruence for the coefficients of the characteristic polynomial of $\overline{\alpha}$, i.e. of the $L$-function $L(P,\chi_d^\kappa;T)$.

\medskip

{\bf Definition 2.3.} {\it Let $<\kappa>_p$ be the subset of $\{1,\dots,d-1\}$ containing all the $\kappa_s$, $0\leq s\leq m-1$, i.e. the set of representatives for the elements of the orbit of multiplication by $p$ in $\Z/d\Z$ containing $\kappa$. From $<\kappa>_p$, we define
$$Y_n(\kappa):=\sum_{s=0}^{\#<\kappa>_p-1} Y_n^{(s)}(\kappa),$$
and the polynomial 
$$\P_{n,\kappa}(\widehat{a}_1,\dots,\widehat{a}_e):=\prod_{s=0}^{\#<\kappa>_p-1}\P_{n,\kappa}^{(s)}(\widehat{a}_1,\dots,\widehat{a}_e).$$}

\medskip 

With these definitions, we can give an expression for the ``principal" part of the minors $\M_n$ of the matrix of $\overline{\alpha}$, i.e. for the coefficients of the $L$-function.

\medskip

{\bf Proposition 2.3.} {\it If $p\geq 2de$, $d\geq 3$, and $\M_n$ denotes the coefficient of degree $n$ of $L(P,\chi_d^\kappa,T)$, we have the following congruence
$$\M_n\equiv \frac{\prod_{s=0}^{m-1}\P_{n,\kappa}^{(s)}(\widehat{a}_1,\dots,\widehat{a}_e)^{\tau^s}}{\left(\prod_{s=0}^{\#<\kappa>_p-1}\prod_{i\notin 
B_n^{(s)}(\kappa)}\lceil\frac{pi-K_s}{e}\rceil!\prod_{i\in 
B_n^{(s)}(\kappa)}\left(\lceil\frac{pi-K_s}{e}\rceil-1\right)!\right)^{\frac{m}{\#<\kappa>_p}}}\pi^{\frac{m}{\#<\kappa>_p}Y_n(\kappa)}$$
modulo $[\pi^{\frac{m}{\#<\kappa>_p}Y_n(\kappa)+1}]$. Moreover, we have $v(\M_n)=\frac{m}{\#<\kappa>_p}Y_n(\kappa)$ if and only if $\P_{n,\kappa}(\widehat{a}_1,\dots,\widehat{a}_e)$ is not zero modulo $p$.}

\medskip

{\it Proof.} From Proposition 2.2, and reasoning as in the proof of Proposition 3.2 in \cite{bf}, we get 
$$\M_n\equiv \frac{\prod_{s=0}^{m-1}\P_{n,\kappa}^{(s)}(\widehat{a_1},\dots,\widehat{a_e})}{\prod_{s=0}^{m-1}\prod_{i\notin 
B_n^{(s)}(\kappa)}\lceil\frac{pi-K_s}{e}\rceil!\prod_{i\in 
B_n^{(s)}(\kappa)}\left(\lceil\frac{pi-K_s}{e}\rceil-1\right)!}\pi^{Y_n}\quad[\pi^{Y_n+1}],$$
with $Y_n=\sum_{s=0}^{m-1} Y_n^{(s)}(\kappa)$. Now remark that since $p^m\equiv1~[d]$, the permutation induced by multiplication by $p$ modulo $d$ has order dividing $m$; thus the length of the cycle containing $\kappa$ divides $m$. Moreover, we get $\kappa_{s+\#<\kappa>_p}=\kappa_s$, and $K_{s+\#<\kappa>_p}=K_s$, $Y_n^{(s+\#<\kappa>_p)}(\kappa)=Y_n^{(s)}(\kappa)$, etc... Consequently $Y_n=\frac{m}{\#<\kappa>_p}Y_n(\kappa)$, and so on for the products of factorials at the denominator. This proves the first assertion. 

\medskip

On the other hand we have, for any $s$, $\P_{n,\kappa}^{(s)}(\widehat{a}_1,\dots,\widehat{a}_e)^{\tau^s}=0$ if and only if $\P_{n,\kappa}^{(s)}(\widehat{a}_1,\dots,\widehat{a}_e)=0$. Thus from the first assertion, we have $v(\M_n)=\frac{m}{\#<\kappa>_p}Y_n(\kappa)$ exactly when $\prod_{s=0}^{m-1}\P_{n,\kappa}^{(s)}(\widehat{a}_1,\dots,\widehat{a}_e)$ is non zero modulo $p$. Finally, we have as above $\P_{n,\kappa}^{(s)}=\P_{n,\kappa}^{(s+\#<\kappa>_p)}$, and consequently $\prod_{s=0}^{m-1}\P_{n,\kappa}^{(s)}=\P_{n,\kappa}^\frac{m}{\#<\kappa>_p}$, and this gives the second assertion since $\P_{n,\kappa}(\widehat{a}_1,\dots,\widehat{a}_e)$ lives in the ring of integers of $\K_m$, an unramified extension of $\Q_p$ in which the ideal $(p)$ is maximal.

\medskip

\section{Newton polygons of the $L$-functions of twisted sums.}

We come to the Newton polygons: Proposition 2.3 above allows us to give a lower bound, the {\it generic Newton polygon}; we show that this lower bound is always above a fixed polygon, we call the {\it Hodge-Stickelberger polygon} of the situation. We show that these bounds are tight, since the Hodge-Stickelberger polygon coincides with the Newton polygon when $p\equiv 1~[de]$. Moreover the generic Newton polygon is attained for every polynomial in a Zariski dense open subset of $\A^{e-1}$, the affine space of monic polynomials of degree $e$ without constant coefficient.

\medskip

{\bf Definition 3.1.} {\it {\it i)} For every $1\leq \kappa \leq d-1$, set $\mu_\kappa:=\frac{\sum_{i\in<\kappa>_p}i}{d\#<\kappa>_p}$. Let $1\leq r\leq d-1$ be the least positive integer congruent to $p$ modulo $d$. We define the {\rm Hodge-Stickelberger polygon} $HS(d,e,r,\kappa)$ as the polygon of horizontal length $e$ and vertices 
$$\left\{(0,0),\left(n,\frac{1}{e}(\frac{n(n-1)}{2}+n\mu_{d-\kappa})\right)_{1\leq n\leq e}\right\}.$$

{\it ii)} We define the {\rm generic Newton polygon} $GNP(d,e,p,\kappa)$ to be the polygon with vertices $$\left\{(0,0),\left(n,\frac{Y_n(\kappa)}{(p-1)\#<\kappa>_p}\right)_{1\leq n \leq e}\right\}.$$
We denote its slopes by  $\lambda_n(\kappa):=\frac{1}{(p-1)\#<\kappa>_p}(Y_{n+1}(\kappa)-Y_{n}(\kappa))$, $0\leq n\leq e-1$.}

\medskip

Let us explain the terminology we employ, and give the connection between the Hodge-Stickelberger polygon and the polygon defined in \cite{as1}, \cite{as2}.

\medskip

{\bf Remark 3.1.} {\it i) The number $\mu_{d-\kappa}$ is well known since Stickelberger. Notations being as above, consider the Gauss sum over $\F_q$
$$\G_{\ma{F}_q}(\Psi_m,\chi^\kappa):=\sum_{x\in \ma{F}_q^\times} \Psi_m(x)\chi^\kappa(x).$$
This is a Weil number, and Stickelberger's theorem ({\it cf.} \cite{bew} Theorem 11.2.1) gives its $q$-adic valuation, which is exactly $\mu_{d-\kappa}$, as can be seen from \cite{bew} Theorem 11.2.7.

\smallskip

ii) Assume $p$ is {\rm semi primitive modulo $d$}, i.e. that $d-1$ is in the subgroup of $\left(\Z/d\Z\right)^\times$ generated by $p$. In this case, it is clear that for any $i$ in $<\kappa>_p$, $d-i$ is also in $<\kappa>_p$, and we get $\mu_\kappa=\frac{1}{2}$. We can rewrite the vertices of the Hodge-Stickelberger polygon in the following simple way: $\left\{(0,0),(n,\frac{n^2}{2e})_{1\leq n\leq e}\right\}$.

\smallskip

iii) Recall (cf. \cite{ro}) that the Hodge polygon $H(e)$ associated to additive character sums and polynomials of degree $e$ is the polygon with vertices $\left\{(0,0),\left(n,\frac{n(n+1)}{2e}\right)_{1\leq n\leq e}\right\}$. We see that the Hodge-Stickelberger polygon defined above differs from $H(e)$ by adding a segment of length one and slope $\frac{\mu_{d-\kappa}}{e}$, and adding $\frac{\mu_{d-\kappa}}{e}$ to each of the slopes.

\medskip

iv) When $\kappa=0$, the polygon with vertices $\left\{(0,0),\left(n,\frac{Y_n(0)}{p-1}\right)_{1\leq n \leq e-1}\right\}$ is exactly the generic Newton polygon in \cite{bf} Definition 4.3.

\medskip

v) The Hodge Stickelberger polygon is the same as the one-variable case of the polygon defined in \cite{as2} Corollary 3.18. Set $\d:=-\frac{(q-1)\kappa}{d}$ in order to comply with the notations there. Then $\d^{(i)}=-\frac{(q-1)\kappa_{m-i}}{d}$, $u_{\bo{d}^{(i)}}(j)=x^{\frac{\bo{d}^{(i)}}{q-1}+j}$, and $w(u_{\bo{d}^{(i)}}(j))=\frac{\frac{\bo{d}^{(i)}}{q-1}+j}{e}=\frac{j}{e}-\frac{\kappa_{m-i}}{de}$ since the weight of $x^r$ is $\frac{r}{e}$ in our case. The $b_j$ defined above Theorem 3.17 in \cite{as2} can be rewritten in the following way 
$$b_j=\frac{1}{m}\left(\sum_{i=0}^{m-1}\frac{j}{e}-\frac{\kappa_{m-i}}{de}\right)=\frac{j-\mu_\kappa}{e}=\frac{j-1+\mu_{d-\kappa}}{e},$$
which gives exactly the $j$-th slope of the Hodge-Stickelberger polygon.}

\medskip

We begin by a lemma, in order to show that the generic Newton polygon as defined above is convex.

\medskip

{\bf Lemma 3.1.} {\it Assume $p\geq 2de$; then for any $0\leq n\leq e-2$, we have $\lambda_n(\kappa)<\lambda_{n+1}(\kappa)$.}

\medskip

{\it Proof.} We have to show the following inequality, for any $0\leq n\leq e-2$: $2Y_{n+1}(\kappa)<Y_{n}(\kappa)+Y_{n+2}(\kappa)$. It is sufficient to show the inequality $2Y_{n+1}^{(s)}(\kappa)<Y_{n}^{(s)}(\kappa)+Y_{n+2}^{(s)}(\kappa)$ for any $s$. We use the expression in Lemma 2.2; first remark that $\#B_n^{(s)}\leq \#B_{n+1}^{(s)}\leq \#B_n^{(s)}+2$. From this we get
$Y_{n+2}^{(s)}(\kappa)\geq Y_{n+1}^{(s)}+\lceil\frac{p(n+2)-K_s}{e}\rceil-2$, and we are reduced to show that $Y_{n+1}^{(s)}(\kappa)<Y_{n}^{(s)}+\lceil\frac{p(n+2)-K_s}{e}\rceil-2$. Now since $Y_{n+1}^{(s)}(\kappa)\leq Y_{n}^{(s)}+\lceil\frac{p(n+1)-K_s}{e}\rceil$, we have to show $\lceil\frac{p(n+2)-K_s}{e}\rceil-\lceil\frac{p(n+1)-K_s}{e}\rceil\geq 2$, and this is verified as long as $p\geq 3e$.

\medskip

In the next proposition, we show that the Hodge-Stickelberger polygon is a lower bound for the Newton polygons of $L$-functions associated to twisted exponential sums as above, and that for $p\equiv 1 ~[de]$, this lower bound is exactly the Newton polygon.

\medskip

{\bf Proposition 3.1.} {\it Let $p\geq 2de$ be a prime, and $1\leq r\leq d-1$ be the least positive integer congruent to $p$ modulo $d$. 
\begin{itemize}
	\item[i)]For any polynomial $P$ of degree $e$ over $\F_q$, the Newton polygon of $L(P,\chi_d^\kappa;T)$ with respect to the $q$-adic valuation lies above the Hodge-Stickelberger polygon $HS(d,e,r,\kappa)$;
	\item[ii)] If $p\equiv 1 ~[de]$, then for any polynomial $P$ of degree $e$ over $\F_q$, the Newton polygon of $L(P,\chi_d^\kappa;T)$ with respect to the $q$-adic valuation is exactly the Hodge-Stickelberger polygon $HS(d,e,r,\kappa)$. 
\end{itemize}}

\medskip

{\it Proof.} From Proposition 2.3, we see that the Newton polygon of $L(P,\chi_d^\kappa,T)$ with respect to the $\pi$-adic valuation lies above the polygon with vertices 
$$\left\{(0,0),\left(n,\frac{m}{\#<\kappa>_p}Y_n(\kappa)\right)_{1\leq n \leq e}\right\}.$$ 
Now we have, for any $0\leq s\leq m-1$, $Y_n^{(s)}(\kappa)\geq \sum_{i=1}^n \frac{pk-\sigma(k)-K_s}{e}=(p-1)\frac{n(n+1)}{2e}-n\frac{K_s}{e}$. Thus we have $\frac{m}{\#<\kappa>_p}Y_n(\kappa)\geq m(p-1)\frac{n(n+1)}{2e}-n\frac{\sum_{s=0}^{m-1}K_s}{e}$. From the definition of $K_s$, we get $\sum_{s=0}^{m-1}K_s=\sum_{s=0}^{m-1}\frac{p\kappa_{s+1}-\kappa_s}{d}=\frac{p-1}{d}\sum_{s=0}^{m-1}\kappa_s$ since $\kappa_m=\kappa_0=\kappa$. The last sum equals $\frac{m}{\#<\kappa>_p}\sum_{i\in<\kappa>_p}i=md\mu_\kappa$, and we get $\frac{m}{\#<\kappa>_p}Y_n(\kappa)\geq m(p-1)\left(\frac{n(n+1)}{2e}-\frac{n\mu_\kappa}{e}\right)$. In other words, we have $v_q(\M_n)\geq \frac{n(n+1)}{2e}-\frac{n\mu_\kappa}{e}=\frac{n(n-1)}{2e}+n\frac{1-\mu_\kappa}{e}$. Finally, $1-\mu_\kappa=1-\frac{\sum_{i\in<\kappa>_p}i}{d\#<\kappa>_p}=\frac{\sum_{i\in<\kappa>_p}d-i}{d\#<\kappa>_p}=\mu_{d-\kappa}$, and this shows part {\it i)} of the proposition.

\medskip

Assume $p\equiv 1 ~[de]$; in this case for any $0\leq s\leq m-1$, we have $\kappa_s=\kappa$, $K_s=\frac{p-1}{d}\kappa$, and for any $1\leq \kappa\leq d-1$, $\mu_\kappa=\frac{\kappa}{d}$, since $<\kappa>_p=\{\kappa\}$. Thus $\lceil\frac{pi-\sigma(i)-K_s}{e}\rceil=\frac{p-1}{e}i-\frac{p-1}{de}\kappa+\lceil\frac{i-\sigma(i)}{e}\rceil$. Summing over $i\in\{1,\dots,n\}$, we see that the sum is minimal exactly when $\sigma(i)\geq i$ for any $i$. But this is possible only when $\sigma$ is the identity of $S_n$. Consequently $\Sigma_n^{(s)}(\kappa)=\{Id\}$, $Y_n^{(s)}(\kappa)=\frac{p-1}{e}\frac{n(n+1)}{2}-\frac{p-1}{de}n\kappa=\frac{p-1}{e}\frac{n(n-1)}{2}-\frac{p-1}{de}n(d-\kappa)$, $\P_{n,\kappa}^{(s)}(X_1,\dots,X_e)=X_e^{Y_n^{(s)}(\kappa)}$, and $\P_{n,\kappa}(X_1,\dots,X_e)=X_e^{Y_n(\kappa)}$. From this and proposition 2.3 we deduce that for any $P$ of degree $e$, the Newton polygon of $L(P,\chi_d^\kappa;T)$ with respect to the $q$-adic valuation is exactly the Hodge-Stickelberger polygon defined above. This proves part {\it ii)}.

\medskip

We see that the situation is rather simple when $p\equiv 1~[de]$: in this case the generic Newton and Hodge-Stickelberger polygons coincide, and the generic Newton polygon is the Newton polygon for any $P$. In the general case, however, the Newton polygon can vary much with the polynomial $P$; we now show that there is a certain regularity: ``most of" the polynomials of fixed degree share the same Newton polygon, which is the generic Newton polygon.

\medskip

{\bf Definition 3.2.} {\it Recall that if $\P(X)=\sum_{i=1}^e \widehat{a_i}X^i$, we have set 
$$\P_{n,\kappa}^{(s)}(\widehat{a}_1,\dots,\widehat{a}_e):=\sum_{\sigma\in \Sigma_n^{(s)}} \sgn(\sigma)
\prod_{i=1}^n\left\{\P^{\lceil\frac{pi-\sigma(i)-K_s}{e}\rceil}\right\}_{pi-\sigma(i)-K_s},$$
and $\P_{n,\kappa}(\widehat{a}_1,\dots,\widehat{a}_e):=\prod_{s=0}^{\#<\kappa>_p-1}\P_{n,\kappa}^{(s)}(\widehat{a}_1,\dots,\widehat{a}_e)$.
\begin{itemize}
	\item[i)] We denote by $P_{n,\kappa}\in \F_p[X_1,\dots,X_e]$ the reduction modulo $p$ of $\P_{n,\kappa}$ ;
	\item[ii)] let 
	$$P_{d,e,p,\kappa}:=\prod_{n=1}^e P_{n,\kappa};$$
we call this polynomial the {\rm Hasse polynomial} attached to the orbit $<\kappa>_p$.
\end{itemize}}
  
\medskip

Reasoning as in Lemma 4.1 of \cite{bf}, we see that the polynomials $P_{n,\kappa}^{(s)}$ are always non zero, and that the polynomial $P_{d,e,p,\kappa}$ itself is non zero. Now from proposition 2.3 and the definition above, we have $P_{d,e,p,\kappa}(a_1,\dots,a_{e-1},1)\neq 0$ exactly when the Newton polygon of $L(P,\chi_d^\kappa;T)$, $P(X)=X^e+a_{e-1}X^{e-1}+\dots+a_1X$, with respect to the $\pi$-adic valuation coincides with the polygon with vertices $\{(0,0),(n,\frac{m}{\#<\kappa>_p}Y_n(\kappa)
)_{1\leq n \leq e}\}$. Thus we get the following, where we identify the set of monic degree $e$ polynomials without constant coefficient with affine space of dimension $e-1$

\medskip

{\bf Theorem 3.1.} {\it For any $p\geq 2de$, the Newton polygon of $L(P,\chi_d^\kappa,T)$ with respect to the $q$-adic valuation coincides with the generic Newton polygon $GNP(d,e,p,\kappa)$ exactly when the polynomial $P$ in the Zariski dense open subset of $\A_{\overline{\ma{F}}_p}^{e-1}$ defined as the complementary of the hypersurface $P_{d,e,p,\kappa}(X_1,\dots,X_{e-1},1)=0$.}

\medskip

{\bf Remark.} When $\kappa=0$, i.e. when we have a purely additive character sum, we know (cf. \cite{zhu2}) that there is a Zariski dense open subset $\U$ defined over $\Q$ of the space of polynomials of degree $e$ over $\overline{\Q}$ such that when $p$ tends to infinity, for any $f\in \U$, the polygon $NP_q(f)$ obtained from the $L$-function associated to the reduction of $f$ modulo a prime above $p$ in the field defined by the coefficients of $f$ tends to a fixed polygon, the Hodge-Stickelberger polygon of this situation. 

\medskip

When $\kappa\neq 0$, we can't hope such a limit to exist, since the polygon $GNP(d,e,p,\kappa)$ depends heavily on $<\kappa>_p$, i.e. on the residue of $p$ modulo $d$. However one can show from the results above that there exists a polynomial $\P_{d,e,r,\kappa}\in \Q[X]$ such that for any sufficiently great $p$, $p\equiv r~[d]$, the Hasse polynomial $P_{d,e,p,\kappa}$ is the reduction modulo $p$ of $\P_{d,e,r,\kappa}$. On the other hand, when $p$ tends to infinity in the class of $r$ modulo $d$, the generic Newton polygon $GNP_q(d,e,p,\kappa)$ tends to the Hodge-Stickelberger polygon $HS(d,e,r,\kappa)$; thus we get a weaker, but similar result than in the case of additive exponential sums: there is a Zariski dense open subset $\U_{d,e,r,\kappa}$ defined over $\Q$ of the space of polynomials of degree $e$ over $\overline{\Q}$ such that when $p$ tends to infinity along the class of $r$ modulo $d$, for any $P\in \U_{d,e,r,\kappa}$, the Newton polygon of the $L$-function $L(P,\chi_d^\kappa;T)$ associated to the reduction of $P$ modulo a prime above $p$ in the field defined by the coefficients of $P$ tends to the Hodge-Stickelberger polygon $HS(d,e,r,\kappa)$.

\medskip

\section{Newton polygons for polynomials $P(x^d)$.}

We come to Newton polygons of $L$-functions associated to additive character sums and polynomials of the form $P(x^d)$ over $k=\F_q$. The idea is to express these $L$-functions as products of $L$-functions of the form $L(P,\chi_d^\kappa;T)$ over extensions of $k=\F_q$. Note that in this section all the results are valid without the hypothesis that $d$ divides $q-1$.

\medskip

\subsection{Another expression for the $L$-function.}
We begin by looking at the sums 
$$S_r(P(x^d))=\sum_{x\in k_r} \psi_{rm}(P(x^d))=\sum_{x\in k_r^\times} \psi_{rm}(P(x^d))+1,$$
and expressing them from the sums $S_r(P,\chi_d^\kappa)$ when $\kappa$ varies

\medskip

{\bf Lemma 4.1.} {\it For any $r\geq 1$, set $\delta_r:=\gcd(d,q^r-1)$, and let $\chi_r:=\chi_d^{d/\delta_r}$ be a multiplicative character of $k_r^\times$ of order $\delta_r$. Then we have
$$S_r(P(x^d))=\sum_{\kappa=0}^{\delta_r-1}\sum_{x\in k_r^\times} \psi_{rm}(P(x^d))\chi_r^\kappa(x)+1=\sum_{\kappa=0}^{\delta_r-1}S_r(P,\chi_r^\kappa)+1.$$}

\medskip

{\it Proof.} Write $d=\delta_r\epsilon_r$; since $\delta_r=\gcd(d,q^r-1)$, $\epsilon_r$ is prime to $q^r-1$, and the map $x\mapsto x^{\epsilon_r}$ is a bijection on $k_r^\times$. On the other hand, since $\delta_r|q^r-1$, the kernel of the map $x\mapsto x^{\delta_r}$ is the set of $\delta_r$-th roots of unity, and its image is the set $(k_r^\times)^{\delta_r}$ of $\delta_r$-th powers in $k_r^\times$. Thus we get
$$S_r(P(x^d))=\sum_{x\in k_r^\times} \psi_{rm}(P(x^d))+1=\delta_r\sum_{x\in (k_r^\times)^{\delta_r}} \psi_{rm}(P(x))+1.$$
Now let $\chi_r$ be as in the Lemma; from the orthogonality relations on multiplicative characters, we have $\sum_{\kappa=0}^{\delta_r-1}\chi_r(x^\kappa)=\delta_r$ if $x\in (k_r^\times)^{\delta_r}$, $0$ else. Replacing in the sum above we obtain
$$\begin{array}{rcl}
S_r(P(x^d)) & = & \sum_{x\in k_r^\times}\sum_{\kappa=0}^{\delta_r-1}\chi_r(x^\kappa) \psi_{mr}(P(x))+1 \\
& = & \sum_{\kappa=0}^{\delta_r-1}\sum_{x\in k_r^\times}\chi_r(x^\kappa) \psi_{mr}(P(x))+1 \\
& = & \sum_{\kappa=0}^{\delta_r-1}S(P,\chi_r^\kappa)+1, \\
\end{array}$$
and this what we wanted to show.

\medskip

From the lemma above, we deduce a factorisation of the $L$-function; let us begin by setting some notations.

\medskip

{\bf Definition 4.1.} {\it i) In $\Z/d\Z$, we choose a set of representatives $\{\kappa_{0,q},\dots,\kappa_{t,q}\}$ for the orbits under multiplication by $q$, and we denote by $\Z/d\Z:=\coprod_{i=0}^t <\kappa_{i,q}>_q$ the decomposition of $\Z/d\Z$ in orbits under multiplication by $q$. For any $i\in \{0,\dots,t\}$, let $d_i:=\#<\kappa_{i,q}>_q$. 

\smallskip

ii) For any $i$, we denote by $L(P,\chi_d^{\kappa_{i,q}};T)$ the $L$-function associated to the character sums 
$$\sum_{x\in k_{d_ir}^\times}\chi_d^{\kappa_{i,q}}(x) \psi_{md_ir}(P(x)).$$}

\medskip

Note that any two elements $\kappa,\kappa'$ of $<\kappa_{i,q}>_q$ are in the same orbit under multiplication by $q$; thus they are in the same orbit under multiplication by $p$, $\mu_{\kappa}=\mu_{\kappa'}$, and this number just depends on the orbit, not on the particular choice of a representative.

\medskip

{\bf Proposition 4.1.} {\it Notations being as above, we have
$$L(P(x^d);T)=\prod_{i=0}^t L(P,\chi_d^{\kappa_{i,q}};T^{d_i}).$$}

\medskip

{\it Proof.} We claim that for any $\kappa,r$ we have $S(P,\chi_r^\kappa)=S(P,\chi_r^{q\kappa})$. Actually we have 
$$S(P,\chi_r^{\kappa})=\sum_{x\in k_r^\times}\chi_r(x^{q\kappa}) \psi_{mr}(P(x^q))$$
since $x\mapsto x^q$ is an automorphism of $k_r$; on the other hand, $P(x^q)=P(x)^q$ since $P\in \F_q[X]$, and this element of $k_r$ has the same trace than $P(x)$; thus $\psi_{mr}(P(x^q))=\psi_{mr}(P(x))$ for any $x$ in $k_r$, proving the claim.

\medskip

From this claim and the result of the lemma above, we can rewrite the sums $S_r(P(x^d))$ using the orbits $<\kappa_{i,q}>_q$. Remark that $d_i|r$ if and only if $q^r\kappa_{i,q}\equiv \kappa_{i,q}~[d]$, i.e. when $d$ divides $(q^r-1)\kappa_{i,q}$, or equivalently when $\frac{d}{\delta_r}|\kappa_{i,q}$. Thus the sum in Lemma 4.1 can be rewritten
$$S_r(P(x^d))=\sum_{\kappa_{i,q},~d_i|r}d_iS_r(P,\chi_d^{\kappa_{i,q}})+1=\sum_{\kappa_{i,q},~d_i|r,~i\neq 0}d_iS_r(P,\chi_d^{\kappa_{i,q}})+S_r(P),$$
and plugging this equality into the formula defining the $L$-function from the character sums gives the desired result.

\medskip

{\bf Remark.} We give below a cohomological interpretation of the equality in the Proposition above. Another way of expressing the $L$-function $L(P(x^d);T)$ is to see it as the characteristic polynomial of the map $\overline{\gamma}$ in the commutative diagram below ({\it cf.} \cite{bf} section 1.2)
$$\xymatrix{
 0 \ar[r]& \ar[d]_{q\gamma} \H^\dagger(A) \ar[r]^{\Delta} & \ar[d]_{\gamma} 
\H^\dagger(A) \ar[r] & \ar[d]_{\overline{\gamma}} \H^\dagger(A)/\Delta\H^\dagger(A) \ar[r] & 0\\
0 \ar[r]& \H^\dagger(A) \ar[r]^{\Delta} & \H^\dagger(A) \ar[r] & \H^\dagger(A)/\Delta\H^\dagger(A) \ar[r] & 0,\\
}$$
where $\Delta:=H(X^d)^{-1}\circ X\frac{d}{dX}\circ H(X^d)$, $\gamma:=H(X^d)^{-1} \circ \psi_q\circ H(X^d)$, and $H(X)=\exp(-\pi\P(X))$. Consider the following decomposition of $\H^\dagger(A)$
$$\H^\dagger(A)=\bigoplus_{\kappa=0}^{d-1} E_\kappa,\qquad E_\kappa=\Vect\{X^{\kappa+dk}\}_{k\in \ma{Z}}.$$
The vector spaces $E_\kappa$ are stable under the action of $\Delta$, and we deduce the following decomposition $\H^\dagger(A)/\Delta\H^\dagger(A)=\bigoplus_{i=0}^{d-1} E_\kappa/\Delta E_\kappa$. On the other hand, the map $\gamma$ sends $E_\kappa$ to $E_{\kappa'}$, and induces $\overline{\gamma}_\kappa$ from $E_\kappa/\Delta E_\kappa$ to $E_{\kappa'}/\Delta E_{\kappa'}$, where $\kappa'$ is the residue of $q\kappa$ modulo $d$. Moreover it can be shown that if $\ell:=\#<d-\kappa>_q$, the map $\overline{\gamma}_{d-\kappa}^\ell$ from $E_{d-\kappa}/\Delta E_{d-\kappa}$ to itself is the same (via a basis change) as the map $\overline{\alpha}$ in section 1.2. Finally an easy calculation gives the following link between the characteristic polynomials of $\overline{\gamma}$ and the $\overline{\gamma}_\kappa$
$$\det\left(\I-T\overline{\gamma}\right)=\prod_{i=0}^t \det\left(\I-T^{d_i}\overline{\gamma}_{\kappa_{i,q}}^{d_i}\right),$$
and we find the same equality as above.

\medskip

\subsection{Newton and Hodge-Stickelberger polygons.} As in section 3, we study the behaviour of the Newton polygon of the $L$-function $L(P(x^d);T)$ when $P$ runs over monic polynomials of degree $e$; we define a Hodge-Stickelberger polygon which is a lower bound for any of these polygons. We show that there is a generic Newton polygon and a Hasse polynomial as above; then studying the case $p\equiv 1~[de]$, we find a classical result of Robba.

\medskip

{\bf Definition 4.2.} {\it i) In $\Z/d\Z$, we choose a set of representatives $\{\kappa_{0,p},\dots,\kappa_{u,p}\}$ for the orbits under multiplication by $p$, and we denote by $\Z/d\Z:=\coprod_{i=0}^u <\kappa_{i,p}>_p$ the decomposition of $\Z/d\Z$ in orbits under multiplication by $p$, with $<\kappa_{0,p}>_p=\{0\}$. For any $i\in \{0,\dots,u\}$, let $D_i:=\#<\kappa_{i,p}>_p$. Assume moreover that the orbits $<\kappa_{i,p}>_p$ are arranged in such a way that $\mu_{\kappa_{1,p}}\leq \dots \leq \mu_{\kappa_{u,p}}$, where the $\mu_{\kappa_{i,p}}$ are as in definition 3.1. We define the {\rm Hodge-Stickelberger polygon} $HS(d,e,r)$ to be the polygon with segments of length $D_i$, and slopes $s_{ij}=\frac{j}{e}+\frac{\mu_{\kappa_{i,p}}}{e}$, $0\leq i\leq u$, $0\leq j\leq e-1$, $(i,j)\neq (0,0)$ rearranged in increasing order of their slopes.

\smallskip

 ii) For $i=0$,  and any $1\leq j\leq e-1$, set $\lambda_j(0):=\frac{1}{p-1}(Y_{j}(0)-Y_{j-1}(0))$; for $i\neq 0$, and $0\leq j\leq e-1$, we set $\lambda_j(\kappa_{i,p}):=\frac{1}{p-1}(Y_{j+1}(\kappa_{i,p})-Y_{j}(\kappa_{i,p}))$. Then we define the {\rm generic Newton polygon} $GNP(d,e,p)$ as the polygon with segments of length $D_i$, and slopes $\lambda_j(\kappa_{i,p})$, $0\leq i\leq u$, $0\leq j\leq e-1$, $(i,j)\neq (0,0)$ rearranged in increasing order of their slopes.} 

\medskip

Note that since we assumed $\mu_{\kappa_{1,p}}\leq \dots \leq \mu_{\kappa_{u,p}}$ and all these numbers are less than $1$, we have $s_{ij}\leq s_{i'j'}$ exactly when $(j,i)\prec(j',i')$ in the lexicographic order. Note also that when $p\equiv 1~[d]$, we have $HS(d,e,1)=H(de)$, the Hodge polygon associated to additive character sums and polynomials of degree $de$, since $HS(d,e,r)$ is the polygon of length $de-1$ and segments with slopes $\frac{i}{de}$, $1\leq i\leq de-1$. Finally, remark that once again the Hodge-Stickelberger polygon above depends only on the residue of $p$ modulo $d$, and not on the particular power $q$ of $p$. Now we have results very similar to the ones in section 3.

\medskip

{\bf Proposition 4.2.} {\it Let $p\geq 2de$ be a prime, and $1\leq r\leq d-1$ be the least positive integer congruent to $p$ modulo $d$. 
\begin{itemize}
	\item[i)]For any polynomial $P$ of degree $e$ over $\F_q$, the Newton polygon of $L(P(x^d);T)$ with respect to the $q$-adic valuation lies above the Hodge-Stickelberger polygon $HS(d,e,r)$;
	\item[ii)] If $p\equiv 1 ~[de]$, then for any polynomial $P$ of degree $e$ over $\F_q$, the Newton polygon of $L(P(x^d);T)$ with respect to the $q$-adic valuation is exactly the Hodge polygon $HS(d,e,1)=H(de)$. 
\end{itemize}}

\medskip

{\it Proof.} When $\kappa=0$, it is well known that the Newton polygon of $L(P(x);T)$ is above the Hodge polygon with segments of length $1$ and slopes $\frac{j}{e}$, $1\leq j\leq e-1$. Now fix a $\kappa_{i,q}\neq 0$; remark that the Newton polygon of $L(P,\chi_d^{\kappa_{i,q}};T)$ with respect to the $q^{d_i}$-adic valuation and of $L(P,\chi_d^{\kappa_{i,q}};T^{d_i})$ with respect to the $q$-adic valuation have the same slopes, but that the length of the segments are multiplied by $d_i$ from the first to the second. Thus this Newton polygon is above the polygon with segments of length $d_i$ and slopes $\frac{j}{e}+\frac{\mu_{\kappa_{i,q}}}{e}$, $0\leq j\leq e-1$.

\medskip

Now remark that we have $<\kappa_{i,q}>_q\subset <\kappa_{i,q}>_p$. Since the permutation induced by multiplication by $q$ is a power of the one induced by multiplication by $p$, we must have $d_i|D_i$, and for any other orbit $<\kappa_{i',q}>_q\subset <\kappa_{i,q}>_p$, the Newton polygons with respect to the $q$-adic valuation of $L(P,\chi_d^{\kappa_{i,q}};T^{d_i})$ and $L(P,\chi_d^{\kappa_{i',q}};T^{d_i})$ have the same lower bound from Proposition 3.1, and the product 
$$L_{\kappa_{i,p}}(T):=\prod_{<\kappa_{j,q}>_q\subset <\kappa_{i,q}>_p}L(P,\chi_d^{\kappa_{j,q}};T^{d_j})$$
has its Newton polygon above the polygon with segments of length $D_i$ and slopes $\frac{j}{e}+\frac{\mu_{\kappa_{i,p}}}{e}$, $0\leq j\leq e-1$, i.e. the image of the polygon $HS(d,e,r,\kappa_i)$ by the similitude of center $(0,0)$ and scale $D_i$.

\medskip

We can rewrite Proposition 4.1 in the following way
$$L(P(x^d);T)=\prod_{i=0}^u L_{\kappa_{i,p}}(T),$$
and we see that the Newton polygon of $L(P(x^d);T)$ is the concatenation of the Newton polygons of the polynomials $L_{\kappa_{i,p}}(T)$, $0\leq i \leq u$ where the segments are reordered with non decreasing slopes. Now the Hodge-Stickelberger polygon $HS(d,e,r)$ is the joint of the images of the polygons $HS(d,e,r,\kappa_{i,p})$ by the similitude of center $(0,0)$ and scale $D_i$, for $0\leq i\leq u$. Since each of these polygons is above the corresponding Hodge-Stickelberger polygon $HS(d,e,r,\kappa_{i,p})$, this ends the proof of Proposition 4.2.

\medskip

{\bf Definition 4.3.} {\it The {\rm Hasse polynomial for polynomials $P(x^d)$ over $\overline{\F}_p$}, $\deg(P)=e$ is
$$ P_{d,e,p}:=P_{d,p}\prod_{i=1}^u P_{d,e,p,\kappa_{i,p}},$$
where $P_{d,p}$ is the polynomial defined in \cite{bf}, Definition 4.1, and the product is taken over the orbits different from $\{0\}$ under multiplication by $p$ in $\Z/d\Z$.}
  
\medskip

As in the case of twisted exponential sums, the polynomial $P_{d,e,p}$ is non zero, and tells us exactly which polynomials have the Newton polygon of their $L$-function attaining the generic Newton polygon. Once more, we identify the set of monic degree $e$ polynomials without constant coefficient with affine space of dimension $e-1$.

\medskip

{\bf Theorem 4.1.} {\it For any $p\geq 2de$, the Newton polygon of $L(P(x^d);T)$ with respect to the $q$-adic valuation coincides with the generic Newton polygon $GNP(d,e,p)$ exactly when the polynomial $P$ is in the Zariski dense open subset $U_{d,e,p}$ of $\A_{\overline{\ma{F}}_p}^{e-1}$ defined as the complementary of the hypersurface $P_{d,e,p}(X_1,\dots,X_{e-1},1)=0$.}

\medskip

{\it Proof.} When $\kappa=0$, we see from \cite{bf} Theorem 4.1 that the Newton polygon of $L(P(x);T)$ coincides with the generic Newton polygon $GNP(e)$ exactly when $P$ is in the complementary of the hypersurface of equation $P_{d,p}=0$, and $GNP(e)$ has slopes $\lambda_1(0),\dots,\lambda_{e-1}(0)$.

\medskip

When $\kappa\neq 0$, we see from Definition 2.3 that the polynomial $P_{d,e,p,\kappa}$ just depends on the orbit $<\kappa>_p$, and not on the choice of $\kappa$ inside the orbit; from Theorem 3.1, we see that for any $\kappa_{i,q}$ such that $<\kappa_{i,q}>_q\subset<\kappa>_p$, the Newton polygon with respect to the $q$-valuation of $L(P,\chi_d^{\kappa_{i,q}};T^{d_i})$ is the polygon with vertices 
$$\left\{(0,0),\left(nd_i,d_i\frac{Y_n(\kappa_{i,q})}{(p-1)D_i}\right)_{1\leq n\leq e}\right\}$$
exactly when $P_{d,e,p,\kappa}(a_1,\dots,a_{e-1},1)\neq 0$. In this case, the Newton polygon of the polynomial $L_{\kappa_{i,p}}(T)$ defined in the proof of Proposition 4.2 has vertices as above, where $D_i$ replaces $d_i$. Reasoning as in the end of the proof of Proposition 4.2, we get the result.

\end{document}